\documentclass{article}
\usepackage{graphicx} %

\title{
A Kaczmarz-Inspired Method for Orthogonalization
}

\author{Isabel Detherage \\ UC Berkeley \and Rikhav Shah\footnote{Supported by NSF CCF-2420130} \\ UC Berkeley}
\date{\today}
\usepackage{graphicx} %
\usepackage{amsmath,amssymb,amsthm}
\usepackage{algpseudocode}
\usepackage{algorithm}
\usepackage{fullpage}
\usepackage{hyperref}
\usepackage{xcolor}
\usepackage{mleftright}
\usepackage{enumitem}

\newcommand{\C}{{\mathbb{C}}}

\newcommand{\Z}{{\mathbb{Z}}}

\newcommand{\eps}{\varepsilon}

\newcommand{\inr}[2]{\left\langle#1,#2\right\rangle}

\newcommand{\wt}{\widetilde}

\newcommand{\qand}{\quad\And\quad}

\newcommand{\abs}[1]{\mleft|#1\mright|}
\newcommand{\magn}[1]{\left\|#1\right\|}

\newcommand{\pare}[1]{\mleft(#1\mright)}

\newcommand{\set}[1]{{\left\{{#1}\right\}}}
\newcommand{\bmat}[1]{\begin{bmatrix}#1\end{bmatrix}}
\newcommand{\alg}[1]{\textnormal{\texttt{#1}}}

\newcommand{\ceil}[1]{\mleft\lceil#1\mright\rceil}
\newcommand{\floor}[1]{\mleft\lfloor#1\mright\rfloor}

\newcommand{\spliteq}[2]{\begin{equation}#1\begin{split}#2\end{split}\end{equation}}

\DeclareMathOperator*{\E}{\mathbb{E}}

\DeclareMathOperator*{\poly}{poly}
\DeclareMathOperator*{\col}{Col}

\DeclareMathOperator*{\tr}{tr}

\newtheorem{theorem}{Theorem}[section]
\newtheorem{lemma}[theorem]{Lemma}
\newtheorem{proposition}[theorem]{Proposition}
\newtheorem{corollary}[theorem]{Corollary}
\newtheorem{fact}[theorem]{Fact}

\theoremstyle{definition}

\newtheorem{remark}{Remark}

\algdef{SE}[REPEATN]{RepeatN}{End}[1]{\algorithmicrepeat\ #1 \textbf{times}}{\algorithmicend}
\begin{document}
\maketitle

\begin{abstract}
This paper asks if the following iterative procedure approximately orthogonalizes a set of $n$ linearly independent unit vectors while preserving their span: in each iteration, access a random pair of vectors and replace one with the component perpendicular to the other, renormalized to be a unit vector.

We provide a positive answer: any given set of starting vectors converges almost surely to an orthonormal basis of their span. We specifically argue that the $n$-volume of the parallelepiped generated by the vectors approaches 1 (i.e. the parallelepiped approaches a hypercube).
If $A$ is the matrix formed by taking these vectors as columns, this volume is simply $\det(\abs A)$ where $\abs A=(A^*A)^{1/2}$. We show that \(O(n^2\log(1/\pare{\det(\abs A)\varepsilon}))\) iterations suffice to bring ${\det(\abs A)}$ above $1-\eps$ with constant probability.

\end{abstract}

\section{Introduction}
We consider a simple procedure for approximately orthogonalizing of a collection of $n$ linearly independent unit vectors, thought of as columns of a matrix $A\in\C^{d\times n},d\ge n$. One iteration of this procedure is as follows: sample two columns, and replace one with its component perpendicular to the other, renormalized to be a unit vector.
The only fixed points of this operation are matrices with orthonormal columns. We ask the following questions: \textit{does this procedure converge to such matrices? If so, at what rate?}

We give a positive answer to the first question and a bound on the rate. Our approach is to consider the evolution of the $n$-volume of the parallelepiped generated by the vectors as the orthogonalizations are performed. This volume is simply $\det(\abs A)$ where $\abs A=(A^*A)^{1/2}$. We show the $n$-volumne is monotone even for adversarially chosen pairs of columns at each iteration, and logarithm of the $n$-volume increases non-trivially in expectation when the pair of columns is selected randomly.

This procedure bears some resemblance to the Kaczmarz method for solving linear systems. If one is solving $A^* x=0$ using Kaczmarz, a row of $A^*$ is sampled and $x$ projected to be orthogonal to the sampled row; this is then repeated many times. In our case, $x$ is itself one of the columns of $A$, say $x=a_k$ and one applies one step of Kaczmarz to the system $(A^* x)_j=0$ for $j\neq k$ (and then renormalizes).

We have two main results concerning ${\det(\abs{A_t})}$ where $A_t$ is the matrix obtained by applying $t$ iterations to $A$. Proposition \ref{a0} provides a bound on $\E\log{\det(\abs{A_t})}$, and Corollary \ref{a1} gives a tail bound for values of $\det(\abs{A_t})$ close to 1. The upshot is that $t=\Theta(n^2\log(1/\pare{\det(\abs A)\eps}))$ iterations suffice to produce ${\det(\abs{A_t})}\ge1-\eps$ with constant probability. This in turn implies, for example, that the nonzero singular values of $A_t$ are contained in $[1-O(\sqrt\eps),1+O(\sqrt\eps)]$, and further that there exists an orthonormal basis $B=\bmat{b_1&\cdots&b_n}$ for $\col(A)$ with
\(\magn{A_t-B}_F\le O(\sqrt\eps)\) (see Facts \ref{a2}, \ref{a3}).

\vspace{5mm}
\noindent\textbf{Concurrent and independent work:} This project was prompted by the very interesting related works of Stefan Steinerberger which inspired us to explore variants of Kaczmarz, particularly those which modify the matrix as the solver runs \cite{b0,b1}. He and the present authors independently conceived of the particular variant analyzed in this paper. %
He recently announced some progress on this problem \cite{b2}. First, he gives conditions on $x$ under which $\magn{Ax}$ increases in expectation. Second, he finds a heuristic argument and numerical evidence that the rate of convergence should be $\kappa(A_t)\sim\exp(-t/n^2)\kappa(A)$ where $A_t$ is the matrix after $t$ orthogonalizations. Our work is independent of his.

\subsection{Technical overview}
We first precisely define the procedure we analyze.
Throughout this paper, let $A=\bmat{a_1&\cdots &a_n}\in\C^{d\times n},d\ge n$ be the decomposition of $A$ into columns.
We assume the columns are unit length and linearly independent.
We define the operation
\spliteq{\label{a4}}{
\alg{orth}&:\C^{d\times n}\times\set{(i,j)\in[n]^2:i\neq j}\to\C^{d\times n}
\\
k\text{th column of }\alg{orth}(A,i,j)&=\begin{cases}
    \hfil a_k & k\neq i\\
    \dfrac{a_i-\inr{a_i}{a_j}a_j}{\magn{a_i-\inr{a_i}{a_j}a_j}} & k=i.
\end{cases}}
At each timestep $t$, our procedure samples $(i_t,j_t)$ uniformly at random from ordered pairs of distinct indices and updates $A_{t+1}\gets \alg{orth}(A_t,i_t,j_t)$.
We study the evolution of the non-negative potential function
\spliteq{\label{a5}}{
\Phi(A)=-\log{\det(\abs A)}\quad\text{where}\quad\abs A=(A^*A)^{1/2}.}
Note that since $A$ has unit-length columns, $\Phi(A)=0$ if and only if the columns of $A$ are orthogonal. We show $\Phi(A_t)$ is monotonically decreasing no matter the selection of $(i_t,j_t)$ (see Lemma \ref{a6}). Over a random selection of $(i_t,j_t)$, it strictly decreases in expectation unless $\Phi(A_t)=0$ already (see Lemma \ref{a7}). Our bound roughly shows that when $\Phi(A_t)\ge 1$, one should expect $\Phi(A_{t+1})\le\Phi(A_t)-O(1/n^2)$, i.e. steady additive progress is made, and when $\Phi(A_t)\le1$, one should expect $\Phi(A_{t+1})\le(1-1/n^2)\Phi(A_t)$, i.e. steady multiplicative progress is made.

\subsection{Algorithmic applications}
The update $A_{t+1}=\alg{orth}(A_t,i_t,j_t)$ can be expressed as
\[A_{t+1}=A_tC_t\]
where $C_t$ is the appropriate column operation, taking just $O(d)$ time to perform.
The corresponding decomposition of $A$ as
\[A=A_\tau C\quad\text{where}\quad C=C_{\tau-1}^{-1}C_{\tau-2}^{-1}\cdots C_0^{-1}\]
is analogous to an approximate QR decomposition: the columns of $A_\tau$ form a well-conditioned basis for $\col(A)$ when $\tau$ is sufficiently large.
Our bounds show $\tau=O(n^2\log(1/\det(\abs A)))$ iterations suffice to achieve a constant condition number; when $\det(\abs A)=\poly(n)$, this gives a total arithmetic cost of $\wt O(dn^2)$, comparable to the arithmetic cost of Gram-Schmidt.

If one requires the decomposition $A=A_\tau C$ have $C$ be upper triangular, one can simply sample the update pairs $(i_t,j_t)$ subject to $i_t>j_t$. This does not change any of the analysis as $\E\abs{\inr{a_i}{a_j}}^2$ is the same. In this case, $A_\tau C$ would in fact be converging to the actual $A=QR$ decomposition.

There are a couple possible advantages to this method over Gram-Schmidt. First, one can end the process early when $A_\tau$ has a ``reasonable'' condition number depending on one's requirements; in particular, it could make a good preconditioner. The memory access pattern is very similar to that of the Kaczmarz method for solving the system $A^*x=b$; one could imagine updating the rows of $A^*$ (and corresponding entries of $b$) to improve their conditioning at the same time the Kaczmarz solver runs.
The idea is that the condition number gets smaller over time, so perhaps the increase in the rate of convergence of Kaczmarz justifies the additional computational expense. For a characterization of the rate of convergence of Kaczmarz in terms of the condition number, see \cite{b3}.
A second possible advantage is that it's highly parallelizable. Section \ref{a8} suggests three concrete ways to perform up to $O(n)$ applications of $\alg{orth}$ in parallel.

\section{Results}
Throughout the paper, $A,A',A_t\in\C^{d\times n}$ all denote matrices with unit length and linearly independent columns.
\subsection{Analysis of a single step}
We first claim $\det(\abs A)$ is monotone in applications of $\alg{orth}$.
\begin{lemma}\label{a6}
    For any $A$ and distinct indices $i,j$, one has
    \[{\det(\abs{A'})}=\frac{{\det(\abs A)}}{\sqrt{1-\abs{\inr{a_i}{a_j}}^2}}\]
    where $A'=\alg{orth}(A,i,j)$.
\end{lemma}
\begin{proof}
$A'$ can be expressed as the composition of two elementary column operations: first, $E_1$ replaces $a_i$ with the component $a_i-\inr{a_i}{a_j}a_j$ so $\det(E_1)=1$; second, $E_2$ renormalizes the $i$th column to be a unit vector, so
\[\det(E_2)=\frac1{\magn{a_i-\inr{a_i}{a_j}a_j}}=\frac1{\sqrt{1-\abs{\inr{a_i}{a_j}}^2}}.\]
Then $A'=AE_1E_2$ so
\[\det(\abs{A'})=\det(E_2^*E_1^*A^*AE_1E_2)^{1/2}=\abs{\det(E_2)}\det(\abs A)\]
as required.
\end{proof}

From the above lemma, we see that the change in the determinant depends only on $\abs{\inr{a_i}{a_j}}^2$. The next lemma relates the average value of this quantity (over a random selection of $i,j$) to the determinant.

\begin{lemma}\label{a9}
For any $A$,
\[\frac1{n(n-1)}\sum_{i\neq j}\abs{\inr{a_i}{a_j}}^2\ge\frac1{(n-1)^2}(1-\det(\abs A)^2).\]
\end{lemma}
\begin{proof}
Note that the left hand side is simply $\frac1{n(n-1)}\magn{A^*A-I}_F^2$.
    Also note that $\magn A_F^2=n$ since $A$ has unit length columns.
    Consider the optimization program $\min\magn{A^*A-I}_F^2$ subject to $\magn A_F^2=n$ and $\det(\abs A)^2=P$. Using the $\magn A_F^2=n$ constraint, the objective can be rewritten as
    \[\magn{A^*A-I}_F^2
    =\tr(A^*AA^*A)-2\tr(A^*A)+n
    =\magn{A^*A}_F^2-n.\]
    Both the constraints and objective function of this program depend only on the singular values of $A$, so we may rewrite the program in terms of them.
    Let $x_1,\ldots,x_n$ be the squares of the singular values of $A$ so that the program becomes
    \[\min\sum_{j\in[n]}x_j^2-n\]
    subject to
    \[x_j\ge0\qand\sum_{j\in[n]}x_j=n\qand \prod_{j\in[n]}x_j=P.\]
By Lemma \ref{a10}, we may set $a=x_2=\cdots=x_n\ge1$ and $b=x_1\le1$. The resulting program of $a,b$ is
\[\min g(a,b)=(n-1)a^2+b^2-n\]
subject to
\spliteq{\label{a11}}{
0\le b\le1\le a\qand (n-1)a+b=n\qand a^{n-1}b=P.}
The equality constraints together imply that the feasible values of $a$ must be solutions to\[f(a)=a^n-\frac n{n-1}a^{n-1}=-\frac P{n-1}.\]
For every $P\in[0,1]$, there is only one solution to the right of 1, denote it $a(P)$.
Since $f(\cdot)$ is monotone to the right of 1, $a(\cdot)$ monotone as well with $a(1)=1$ and $a(0)=\frac n{n-1}$. Let
\[\wt f(a)=(n-1)(a-1)^2-\frac1{n-1}\] 
and denote the largest solution to $\wt f(a)=-\frac P{n-1}$ by $\wt a(P)$.
Observe that $f(a)\le\wt f(a)$ on the interval $[1,n/(n-1)]$, so $\wt a(P)\le a(P)$.
The value of the objective function for $(n-1)a+b=n$ is
\spliteq{}{
g(a,b)
  &=-n+(n-1)a^2+b^2
\\&=-n+(n-1)a^2+(n-(n-1)a)^2
\\&=n(n-1)(a-1)^2.}
In particular, this expression is monotone for $a\ge1$. Thus the minimum value of $g$ subject to the constraints \eqref{a11} is
\spliteq{}{
\min g(a,b)
&=g(a(P), n-(n-1)a(P))
\\&=n(n-1)(a(P)-1)^2
\\&\ge n(n-1)(\wt a(P)-1)^2
\\&\ge\frac{n}{n-1}(1-P)}
as required, where the last step used that $\wt f(\wt a(P))=-P/(n-1)$.
\end{proof}

With these two lemmas in place, we can compute the expected value of
\[\Phi(\cdot)=-\log{\det(\abs\cdot)}\] in terms of itself after one application of $\alg{orth}$ for random indices.
Define the iterative map\spliteq{}{f(x)=x-\frac{1-\exp(-2x)}{2(n-1)^2}.}
When $x$ is small, one should think of $f$ as roughly $f(x)\approx(1-O(1/n^2))x$. For larger $x$, one should think of $f$ as roughly $f(x)\approx x - O(1/n^2)$.

\begin{lemma}[One step estimate]\label{a7}
    Let $(i,j)$ be a uniform sample from $\set{(i,j)\in[n]^2:i\neq j}$. Then
    \[\Phi(A')\le\Phi(A)\text{ pointwise}\qand\E\Phi(A')\le f(\Phi(A))\]
where $A'=\alg{orth}(A,i,j)$.
\end{lemma}
\begin{proof}
Lemma \ref{a6} states
\spliteq{\label{a12}}{
\Phi( A')=\Phi(A)+\frac12\log\pare{1-\abs{\inr{a_i}{a_j}}^2}.
}
Since $\abs{\inr{a_i}{a_j}}^2\ge0$, the first claim follows immediately.
By Jensen's inequality, when $i,j$ are picked randomly one has
\spliteq{}{
\E\Phi( A')
  &\le\Phi( A)+\frac12\log\pare{1-\E\abs{\inr{a_i}{a_j}}^2}
\\&\le\Phi( A)-\frac12\E\abs{\inr{a_i}{a_j}}^2
\\&=  \Phi( A)-\frac{\magn{A^*A-I}_F^2}{2n(n-1)}.
}
Now apply Lemma \ref{a9},
    \spliteq{}{
\E\Phi( A')
  &\le\Phi(A)-\frac{1-\det(\abs A)^2}{2(n-1)^2}
\\&=\Phi(A)-\frac{1-\exp(-2\Phi(A))}{2(n-1)^2}
\\&=f(\Phi(A)).}
\end{proof}

\begin{remark}\label{a13}
    Lemma \ref{a7} holds if $(i,j)$ is instead a uniform sample from
    either $\set{(i,j)\in[n]^2:i>j}$ or $\set{(i,j)\in[n]^2:i<j}$.
\end{remark}

\subsection{Supermartingale behavior}
\newcommand{\inflection}{\alpha}
\newcommand{\hit}[2]{{\tau^{({#1})}_{{#2}}}}
\newcommand{\hitt}[2]{\tau_{#1\to#2}}
The previous section describes the effect of one application of $\alg{orth}$ to the collection of vectors. This section considers the effect of several iterations. In particular, we're interested in two dual questions: what does this collection of vectors look like after a given number of iterations? How many iterations are required to achieve a particular state?

For this section, define the family of stochastic processes parameterized by the starting matrix $A$,
\spliteq{\label{a14}}{
A_0=A\qand A_{t+1}=\alg{orth}(A_t,i_t,j_t)
}
where $i_t,j_t$ are sampled uniformly at random (subject to $i_t\neq j_t$) at each step.
We can define the corresponding hitting time of achieving a small potential value when starting with a matrix $A$,
\[\tau^{(A)}_\beta=\inf\set{t\in\Z_{\ge0}:\Phi(A_t)<\beta}.\]
Further define the largest expected hitting time among starting matrices with potential bounded by $\alpha$,
\spliteq{\label{a15}}{
\mu_{\alpha\to\beta}=\sup_{X:\Phi(X)\le\alpha}\E\hit X\beta.
}
This is the most expected amount of time this process takes to transform a matrix with potential at most $\alpha$ to a matrix with potential less than $\beta$.

Because applying $f$ does not commute with taking expectations, obtaining a bound for several applications of $\alg{orth}$ by iteratively applying Lemma \ref{a7} is not immediate. 
Fortunately, $f$ can be bounded by a linear function an adequate slope on bounded domains; this observation results in the following expectation bound.
\begin{proposition}[Expectation bound on $\Phi(A_t)$]\label{a0}
\[\E\Phi(A_t)\le\exp\pare{-\frac{t}{(n-1)^2(2\Phi(A)+1)}}\Phi(A).\]
\end{proposition}
\begin{proof}
Denote $\phi_t=\Phi(A_t)$.
The second part of Lemma \ref{a7} implies
\[\E \phi_{t+1}\le \E f\pare{\phi_t}\]
The first part of Lemma \ref{a7} implies $\phi_t\le\phi_0$ for all $t$ with probability 1. Since $f$ is convex, we have for $\phi\le\phi_0$ that \[f(\phi)\le{\frac{\phi}{\phi_0}}\cdot f(\phi_0).\]
In particular, this inductively implies
\spliteq{\label{a16}}{
\E\phi_{t+1}
\le \frac{f(\phi_0)}{\phi_0}\E\phi_t \le \pare{\frac{f(\phi_0)}{\phi_0}}^{t+1}\phi_0.
}
Finally, one may upper bound $f$ via
\spliteq{}{
f(\phi_0)
  &=  \phi_0-\frac{1-\exp(2\phi_0)^{-1}}{2(n-1)^2}
\\&\le\phi_0-\frac{1-\pare{1+2\phi_0}^{-1}}{2(n-1)^2}
\\&=\phi_0\pare{1-\frac1{(n-1)^2(2\phi_0+1)}}
\\&\le\phi_0\exp\pare{-\frac1{(n-1)^2(2\phi_0+1)}}.}
Combined with \eqref{a16}, this is exactly the desired result.
\end{proof}
\begin{remark}
    For $t=o(n^2\Phi(A))$, the bound of Proposition \ref{a0} resembles linear progress in $t$:
    \[\E\Phi(A_t)\le\Phi(A)-O\pare{\frac t{n^2}}.\]
\end{remark}

This establishes control of just the expectation of $\Phi(A_t)$. The next lemma establishes a tail bound in terms of the $\mu_{\alpha\to\beta}$ defined in \eqref{a15}.

\begin{lemma}[Tail bound in terms of $\mu_{\alpha\to\beta}$]\label{a17}
Fix any $\alpha>\beta>0$. For every $A$ with $\Phi(A)=\alpha$,
\[\Pr(\Phi(A_t)\ge\beta)=\Pr(\hit A\beta>t)\le\exp\pare{-\floor{\frac t{e\mu_{\alpha\to\beta}}}}.\]
\end{lemma}
\begin{proof}
Conditioned on the event $\hit A\beta>k$, we have by definition,
\[\hit A\beta=k+\hit{A_k}\beta.\]
Since $\Phi{(A_t)}$ is pointwise monotone, we have that
\(\Phi(A_k)\le\Phi(A)=\alpha.\)
This yields
\spliteq{}{
\E(\hit{A}\beta\,|\,\hit{A}\beta>k)
  &=\E(\E(k+\hit{A_k}\beta\,|\,\hit{A}\beta>k,A_k))
\\&=k+\E(\E(\hit{A_k}\beta\,|\,\hit{A}\beta>k,A_k))
\\&\le k+\sup_{X:\Phi(X)\le\alpha}\E(\hit{A_k}\beta\,|\,\hit{A}\beta>k,A_k=X)
\\&\le k+\sup_{X:\Phi(X)\le\alpha}\E(\hit X\beta)
\\&=k+\mu_{\alpha\to\beta}.
}
This allows us to compute a bound on the conditional tails of $\hit A\beta$ using Markov's inequality. For $c$ which we specify later,
\spliteq{\label{a18}}{
\Pr\pare{\hit A\beta> k+c\,|\,\hit A\beta> k}
  &=\Pr\pare{\hit A\beta-k> c\,|\,\hit A\beta> k}
\\&\le\frac{\E(\hit A\beta-k\,|\,\hit A\beta> k)}c
\\&\le\frac{\mu_{\alpha\to\beta}}c.}
Set $J=\floor{t/c}$ and apply \eqref{a18} for $k=0,c,2c,\ldots (J-1)c$ to obtain
\spliteq{}{
\Pr\pare{ \hit A\beta> t }
\le
\Pr\pare{ \hit A\beta> Jc }
=\prod_{j=1}^J \Pr\pare{\hit A\beta > jc\,|\,\hit A\beta> (j-1)c}
\le(\mu_{\alpha\to\beta}/c)^J.}
Pick $c=e\mu_{\alpha\to\beta}$ for the final result.\end{proof}

To turn this into a concrete tail bound, we need to establish an estimate for $\mu_{\alpha\to\beta}$. We will do so in two steps: Lemma \ref{a19} establishes a crude upper bound and Lemma \ref{a20} exploits the sub-additive nature of $\mu$ to boost this to a stronger bound.
\begin{lemma}[First expectation bound on $\mu_{\alpha\to\beta}$]\label{a19}For $\alpha>\beta>0$,
\[\mu_{\alpha\to\beta}\le\frac{2(n-1)^2\alpha}{1-\exp(-2\beta)}.\]
\end{lemma}
\begin{proof}
Fix any $A$ with $\Phi(A)\le\alpha$.
Denote $\phi_t=\Phi{(A_t)},\tau_\beta=\hit A\beta$ for brevity.
Note that $f(x)-x$ is negative and monotonically decreasing, so
\spliteq{}{
f(x)-x
&\le\begin{cases}f(\beta)-\beta & x\ge\beta\\0&x<\beta\end{cases}.
}
Using Lemma \ref{a7}, this gives for each $t$ the bound on the increments,
\spliteq{}{
\E(\phi_{t+1}-\phi_t)\le\E(f(\phi_t)-\phi_t)
  &\le
  (f(\beta)-\beta)
  \cdot\Pr\pare{\phi_t\ge\beta}
\\&= 
(f(\beta)-\beta)
\cdot\Pr\pare{\tau_\beta>t}.}
This becomes a telescoping sum, so we obtain for all $t$ that
\[\E(\phi_{t})-\phi_0
\le(f(\beta)-\beta)\sum_{j=0}^{t-1}\Pr\pare{\tau_{\beta}>j }.\]
In the limit as $t\to\infty$, the right hand side is simply $(f(\beta)-\beta)\E(\tau_\beta)$, and by Lemma \ref{a0} the left hand side is $-\phi_0$. By construction, $\phi_0\le\alpha$ so rearranging gives $E(\tau_\beta)\le\alpha/(f(\beta)-\beta)$. Taking the supremum over all $A$ with $\Phi(A)\le\alpha$ gives the final result.
\end{proof}

\begin{lemma}[Better expectation bound on $\mu_{\alpha\to\beta}$]\label{a20}
For $\alpha>1>\beta>0$,
\[\mu_{\alpha\to\beta}\le\frac{2(n-1)^2}{1-e^{-2}}(\alpha+e\ceil{\log(1/\beta)}).\]
For $1>\alpha>\beta>0$,
\[\mu_{\alpha\to\beta}\le\frac{2e(n-1)^2}{1-e^{-2}}\ceil{\log(\alpha/\beta)}.\]
\end{lemma}
\begin{proof}
Consider any sequence $\alpha=\alpha_0>\alpha_1>\cdots>\alpha_\ell$ ending at $\alpha_\ell\le\beta$.
Set $k_0=0$, $k_j=\hit{A_{k_{j-1}}}{\alpha_j}$. Then by definition,
\[\hit A{\beta}\le\sum_{j=1}^\ell k_j.\]
Note also by definition that $\Phi(A_{k_{j-1}})<\alpha_{j-1}$. Thus taking expectations and applying Lemma \ref{a19} gives
\[\E\hit A{\beta}
\le\sum_{j=1}^\ell\E k_j
\le\sum_{j=1}^\ell\mu_{\alpha_{j-1}\to\alpha_j}
\le\sum_{j=1}^\ell\frac{2(n-1)^2\alpha_{j-1}}{1-\exp(-2\alpha_j)}.\]
When $\alpha>1>\beta$, pick $\alpha_j=e^{1-j}$ for $j\ge 1$. Further use the approximation $1-\exp(-2x)\ge(1-e^{-2})x$ when $x\le1$. Then the sum becomes
\spliteq{}{
  \frac{2(n-1)^2\alpha}{1-e^{-2}}
  +
  2(n-1)^2\sum_{j=2}^\ell\frac{e^{2-j}}{1-\exp(-2\cdot e^{1-j})}
&\le  
\frac{2(n-1)^2\alpha}{1-e^{-2}}
+
\frac{2(n-1)^2}{1-e^{-2}}\sum_{j=2}^\ell\frac{e^{2-j}}{e^{1-j}}
\\&=
\frac{2(n-1)^2}{1-e^{-2}}(\alpha+e(\ell-1)).}
It suffices to take $\ell=\ceil{\log(1/\beta)}+1$ giving the first result.
If $1\ge\alpha>\beta$, pick $\alpha_j=\alpha_0e^{-j}$. Then the sum becomes
\[
\sum_{j=1}^\ell\frac{2(n-1)^2\alpha_0e^{1-j}}{1-\exp(-2\alpha_0e^{-j})}
\le
\sum_{j=1}^\ell\frac{2(n-1)^2\alpha_0e^{1-j}}{(1-e^{-2})\alpha_0e^{-j}}
=\frac{2e(n-1)^2}{(1-e^{-2})}\ell.\]
It suffices to take $\ell=\ceil{\log(\alpha/\beta)}$.
\end{proof}

\begin{corollary}[Final tail bound for $\Phi{(A_t)}$]\label{a1}
Fix any $A$ and $\beta<1$. Then
\[\Pr(\Phi(A_t)\ge\beta)\le\exp\pare{-\floor{\frac{1-e^{-2}}{2e}\frac t{(n-1)^2(\Phi(A)+e\ceil{\log(1/\beta)})}}}.\]
In particular,
\[t\ge\frac{2e^3}{e^2-1}(n-1)^2(\Phi(A)+e\ceil{2\log(4/\eps)})\ceil{\log(1/\delta)}=\Theta\pare{n^2\log\pare{\frac1{\abs{\det(A)\eps}}}\log\pare{\frac1\delta}}\]
implies
\[\Pr\pare{\Phi(A_t)\ge1+\eps}\le\delta.\]
\end{corollary}
\begin{proof}
    This is an immediate combination of Lemma \ref{a17} and Lemma \ref{a20}.
\end{proof}

\subsection{Relationship to other measures of orthogonality}
The previous subsection produced results concerning $\Phi(A_t) = -\log\det(\abs{A_t})$. Here we show convergence of $\Phi(\cdot)$ implies convergence of a couple other natural measures of orthogonality.

\begin{fact}[Condition number]\label{a2}
\[\Phi(A)\le\beta\implies\log\kappa(A)\le\min\pare{\log(2)+\beta,\sqrt{2\beta}+\frac{\beta^{1.5}}4}.\]
\end{fact}
\begin{proof}
By applying Lagrange multipliers to the optimization problem
\[\max\kappa(A)=\sigma_1/\sigma_n\]
subject to
\[\magn A_F^2=\sigma_1^2+\cdots+\sigma_n^n=n\qand\sigma_1\cdots\sigma_n=\det(\abs A),\]
one finds that the optimum is achieved for
\[
\sigma_1^2=1+\sqrt{1-\det(\abs A)^2}
\qand
\sigma_2=\cdots=\sigma_{n-1}=1
\qand
\sigma_n^2=1-\sqrt{1-\det(\abs A)^2}
\]
which gives a maximum value of
\[
\max\kappa(A)=\sqrt{
\frac{1+\sqrt{1-\det(\abs A)^2}}{1-\sqrt{1-\det(\abs A)^2}}}\le
\sqrt{
\frac{1+\sqrt{1-\exp(-2\beta)}}{1-\sqrt{1-\exp(-2\beta)}}}\]
the log of which is upper bounded both by $\log(2)+\beta$ and $\sqrt{2\beta}+\beta^{1.5}/4$.\end{proof}

\begin{fact}[Distance to unitary]\label{a3}
\spliteq{}{
\Phi(A)\le\beta&\implies
\inf\set{\magn{A-B}_F:B^*B=I,\col(A)=\col(B)}\le\sqrt{2\beta}.
}
\end{fact}
\begin{proof}
By considering the Gram-Schmidt basis, one may assume without loss of generality that $A$
is upper triangular with positive real diagonal entries.
Let $1=\lambda_1,\cdots,\lambda_n$ be those entries. Note by the Cauchy-Binet formula that $\det(A^*A)=(\lambda_1\cdots\lambda_n)^2$.
Let $r_j\in\C^{j-1}$ be the portion of the $j$th column of $A$ strictly above the diagonal.
Take $B=\bmat{e_1&\cdots&e_n}$. Then
\spliteq{}{
\magn{A-B}_F^2
=\sum_{j=2}^n\pare{(\lambda_j-1)^2+\magn{r_j}^2}
=\sum_{j=2}^n\pare{1-2\lambda_j+\lambda_j^2+\magn{r_j}^2}.
}
Note $\lambda_j^2+\magn{r_j}^2=1$ since the columns of $A$ are unit length. Thus
\spliteq{}{
\magn{A-B}_F^2
=2\sum_{j=1}^n\pare{1-\lambda_j}
\le2\sum_{j=1}^n\log\pare{1/\lambda_j}=-2\log(\lambda_1\cdots\lambda_n)=2\Phi(A).}
\end{proof}

\paragraph{Final thoughts and future directions:}
We do not produce lower bounds, and conjecture that $t=O(n^2\log(\kappa(A)/\eps))$ iterations suffice to bring $\det(\abs{A_t})\ge1-\eps$ with constant probability. Since
\spliteq{\label{a21}}{
\kappa(A)=\frac{\sigma_1(A)}{\sigma_n(A)}\ge\frac1{\det(\abs A)^{1/n}},}
our bound implies only that $O(n^3\log(\kappa(A)/\eps))$ iterations suffice, so a factor of $n$ improvement may be available when \eqref{a21} is loose.

\vspace{5mm}\noindent\textbf{Acknowledgments: }We would like to thank Stefan Steinerberger for his inventive works on the Kaczmarz algorithm and engaging discussions which inspired this project.

\bibliographystyle{alpha}
\bibliography{outbib}

\begin{thebibliography}{Ste21b}

\bibitem[Ste21a]{b0}
Stefan Steinerberger.
\newblock Randomized {K}aczmarz converges along small singular vectors.
\newblock {\em SIAM J. Matrix Anal. Appl}, 42:608--615, 2021.

\bibitem[Ste21b]{b1}
Stefan Steinerberger.
\newblock {Surrounding the solution of a Linear System of Equations from all
  sides}.
\newblock {\em Quarterly of Applied Mathematics}, 79(3):419--429, 2021.

\bibitem[Ste24]{b2}
Stefan Steinerberger.
\newblock Kaczmarz {K}ac walk.
\newblock {\em arXiv preprint 2411.06614}, 2024.

\bibitem[SV09]{b3}
Thomas Strohmer and Roman Vershynin.
\newblock {A Randomized Kaczmarz Algorithm with Exponential Convergence}.
\newblock {\em J Fourier Anal Appl}, 15:262--278, 2009.

\end{thebibliography}

\appendix

\section{Extension to allow parallelism}\label{a8}
The process is currently described in terms of modifying a single column at a time. However, we can batch several of these updates together, computing the changes to multiple columns in parallel.
Each batch of updates corresponds to a directed, rooted forest on the vertex set $[n]$, where edges $(i,j)$ are oriented so that $i$ is further from the root.
Then for each edge $(i,j)$, set the new value of $a_i$ to be
\[\frac{a_i-\inr{a_i}{a_j}a_j}{\magn{a_i-\inr{a_i}{a_j}a_j}}.\]
The orientation of the edges ensures that this is well defined.
By listing the edges of the tree in order of furthest to the roots to the closest, say \((i_0,j_0),\cdots,(i_{m-1},j_{m-1}),\)
one can express the effect of a batch update as $A_m$ where
\[A=A_0\qand A_{t+1}=\alg{orth}(A_t,i_t,j_t)\quad\forall t\in[m-1].\]
Note that once a column is changed, it is not used or changed again. Thus all of the updates to each $A_t$ can be expressed in terms of the columns of the inital $A_0$. In particular, Lemma \ref{a6} simply implies
\spliteq{}{
\Phi(A_t)
  &=\Phi(A)+\frac12\sum_{\ell=1}^t\log(1-\abs{\inr{a_{i_\ell}}{a_{j_\ell}}}^2)
}
analogously to \eqref{a12}. If the forest is sampled so that the marginal distribution of each edge pair $(i_t,j_t)$ is uniform, then following the remainder of the proof of Lemma \ref{a7} gives the generalization
\spliteq{\label{a22}}{
\Phi(A_t)\le\Phi(A)\text{ pointwise }\qand\E\Phi(A_t)\le f_t(\Phi(A))
} where
\[f_t(x)=x-t\cdot\frac{1-\exp(-2x)}{2n(n-1)}.\]
One can then define the process by sampling a sequence of independent forests and concatenating the lists of the edges to form the sequence $\set{(i_t,j_t):t\in\Z_{\ge0}}$. Replacing the result of Lemma \ref{a7} with \eqref{a22} allows the subsequent proofs to go through.

There are many choices of how to sample the required forests. We end this section by suggesting four natural distributions.
\begin{enumerate}
    \item Sample a uniformly random forest consisting of a single edge.
    \item Sample a uniformly random perfect matching of the complete graph and flip a coin to orient each edge.
    \item Sample a uniformly random permutation $\pi\colon[n]\to[n]$ and consider the path with edges $(i_t,j_t)=(\pi(t),\pi(t+1))$ for $t\in[n-1]$.
    \item Sample a uniformly random vertex and consider the star graph rooted and centered at that selected vertex.
\end{enumerate}

\section{Proof of claim in Lemma \ref{a9}}
\label{a23}
\begin{lemma}\label{a10}
    If the program
    \[x_j\ge0\qand\sum_{j\in[n]}x_j=S\qand\prod_{j\in[n]}x_j=P\]
is feasible, then the minimizer of the objective
\[\sum_{j\in[n]}x_j^2\]
occurs for
\[x_{\pi(1)}\le x_{\pi(2)}=\cdots=x_{\pi(n)}\]
for some permutation $\pi$.
\end{lemma}
\begin{proof}
Let $y_1,\ldots,y_n$ be a minimizer of the program and renumber the variables so that $y_1\le y_2\le\cdots\le y_n$ without loss of generality.
    If one fixes $x_j=y_j$ for $j=3,\ldots,n-1$, then the triplet $(y_1,y_2,y_n)$ should minimize the restricted program
    \[\min x_1^2+x_2^2+x_n^2\]
    subject to
    \[
    x_1,x_2,x_n\ge 0
    \qand
    x_1+x_2+x_n=y_1+y_2+y_n=:S'
    \qand
    x_1x_2x_n=y_1y_2y_n=:P'.\]
    Since the constraints and objective are homogeneous, we may assume $S'=3$ without loss of generality.
    Using Lagrange multipliers, one can see that the critical values occur when the set $\set{x_1,x_2,x_n}$ has only two distinct values. Let $a$ and $b$ be the value of the majority and minority element respectively, so that the constraints read
    \[2a+b=3\qand a^2b=P\]
    and the objective is
    \[\min g(a,b)=2a^2 + b^2\]
    Combining the constraints shows that the only feasible values of $a$ are the positive roots of
    \[
    f(a)
    =2a^3-3a^2+P'
    =2(a-1)^3+3(a-1)^2+(P'-1)
    .\]
    Then one has
    \[g(a,3-2a)=3+6(a-1)^2 = 3 + 2\cdot\pare{-(P'-1)-2(a-1)^3}\]
    at those roots. Since $g(a,3-2a)$ is monotonically decreasing in $a$, the minimizer occurs when $a$ is the larger root of $f$. Since the minimum of $f(a)$ for $a\ge0$ occurs at $a=1$, this forces $b\le1\le a$. In particular, we must have $b=y_1\le y_2=\cdots=y_n=a$ as required.
\end{proof}

\end{document}